\input amstex.tex
\documentstyle{amsppt}
\input pictex.tex
\magnification=\magstep1
\pagewidth{15.9truecm}
\nologo
\NoRunningHeads
\refstyle{A}
\widestnumber\key{ACLM}
\topmatter
\title
Stringy invariants of normal surfaces
\endtitle
\author
Willem Veys  \bigskip
\endauthor
\address K.U.Leuven, Departement Wiskunde, Celestijnenlaan 200B,
         B--3001 Leuven, Belgium  \endaddress
\email wim.veys\@wis.kuleuven.ac.be  \newline
{}\quad
 http://www.wis.kuleuven.ac.be/algebra/veys.htm
\endemail
\abstract
The stringy Euler number and $E$--function of Batyrev for log terminal singularities can in dimension 2 also be considered for a normal surface singularity with all log discrepancies nonzero in its minimal log resolution. Here we obtain a structure theorem for  resolution graphs with respect to log discrepancies, implying that these stringy invariants can be defined in a natural way, even when some log discrepancies are zero, and more precisely for {\sl all} normal surface singularities which are not 
log canonical.
We also show that the stringy $E$--functions of log terminal surface singularities are polynomials (with rational powers) with nonnegative coefficients, yielding well defined (rationally graded) stringy Hodge numbers.
\endabstract
\endtopmatter
\bigskip
\document
\head
Introduction
\endhead

\bigskip
\noindent
{\bf 0.1.} The invariants stringy Euler number and $E$--function were defined as follows by Batyrev [B1] for complex algebraic varieties $Y$ with at worst log terminal singularities.  Let $\pi : X \rightarrow Y$ be a log resolution of $Y$ and $E_i, i \in
T$, the irreducible components of the exceptional divisor of $\pi$ with log discrepancies $a_i$, i.e. $K_X = \pi^\ast K_Y + \sum_{i \in T} (a_i - 1)E_i$ where $K_\centerdot$ denotes the canonical divisor.  Denote also $E^\circ_I := (\cap_{i \in I} E_i) \setminus ( \cup_{\ell \not\in I} E_\ell)$ for $I \subset T$.  Then the {\it stringy Euler number} of $Y$ is
$$e(Y) := \sum_{I \subset T} \chi (E^\circ_I) \prod_{i \in I} \frac{1}{a_i} \in \Bbb Q \, , \tag $*$ $$
where $\chi(\cdot)$ denotes the topological Euler characteristic.  Recall that log terminality means that all $a_i > 0$.  A finer invariant is the {\it stringy E--function}
$$E(Y) := \sum_{I \subset T} H(E^\circ_I) \prod_{i \in I} \frac{uv-1}{(uv)^{a_i}-1} \tag $**$ $$
of $Y$, where $H(E^\circ_I) \in \Bbb Z[u,v]$ is the Hodge polynomial of $E^\circ_I$, see (1.1); $E(Y)$ is thus a `rational function with rational powers'.  The proof of Batyrev that these invariants do not depend on the choice of a particular resolution uses the idea of motivic integration, initiated by Kontsevich [K] and developed by Denef and Loeser [DL1,DL2].  In [B1] the stringy $E$--function was used to formulate a topological mirror symmetry test for pairs of Calabi--Yau varieties (with at worst Gorenstein canonical singularities).

These invariants were subsequently generalized to kawamata log terminal pairs [B2], and to pairs $(Y,D)$, where $Y$ is any $\Bbb Q$--Gorenstein variety and $D$ a $\Bbb Q$--Cartier divisor on $Y$ whose support contains the locus of log canonical singularities of $Y$ [V3].
\bigskip
\noindent
{\bf 0.2.}  For normal surfaces $Y$ (with at worst log terminal singularities) it is clear that one can define these invariants using the minimal log resolution of $Y$, and check that they are invariant under blowing--ups.  And of course for any normal surface the expressions ($*$) and ($**$) make sense, at least if all $a_i \ne 0$ ! (This was already implicit in [V2, 5.8].)

Singularities for which certainly not all $a_i \ne 0$ are the log canonical singularities which are not log terminal (we call these {\it strictly log canonical}\,).  But also quite `general' singularities can have some $a_i = 0$ in their minimal log resolution.

Here we will introduce a natural extension of the stringy Euler number and stringy $E$--function to normal surfaces, even when some $a_i = 0$; more precisely we will allow $a_i = 0$ if $E_i \cong \Bbb P^1$ and intersects at most twice other components $E_\ell$.  For example if $E_i$ intersects precisely $E_1$ and $E_2$ (with $a_1 \ne 0$ and $a_2 \ne 0$), then we redefine the contribution of $\{ i \}, \{ i,1 \}$ and $\{i , 2 \}$ in ($*$) as
$$\frac{\kappa_i}{a_1a_2},$$
where $-\kappa_i$ is the self--intersection number of $E_i$ in $X$.  These invariants are defined using the minimal log resolution, but in fact one can use any `allowed resolution' as above.

The point is that in this way we are able to generalize $E(Y)$ and $e(Y)$ to {\it any} normal surface without strictly log canonical singularities.  This is a consequence of the following theorem concerning the structure of the resolution graph of surface singularities.
\bigskip
\proclaim
{0.3. Theorem}  Let $P \in Y$ be a normal surface singularity (germ) which is not log canonical.  Let $\pi : X \rightarrow Y$ be the minimal log resolution of $P \in Y$, for which we use the notation of (0.1).  Then $\pi^{-1} P = \cup_{i \in T} E_i$ consists of the connected part ${\Cal N} := \cup_{\Sb i \in T \\ a_i < 0 \endSb} E_i$, to which a finite number of chains are attached as in Figure 1.
Here $E_i \subset {\Cal N}, E_\ell \cong \Bbb P^1$ for $1 \leq \ell \leq r$, $a_1 \geq 0$ and $(a_i <)a_1 < a_2 < \cdots < a_r <1$.

In particular either all $a_\ell > 0$ for $1 \leq \ell \leq r$, or $a_1 = 0$ and $a_\ell > 0$ for $2 \leq \ell \leq r$.
\endproclaim
\vskip .5truecm
\centerline{
\beginpicture
\setcoordinatesystem units <.5truecm,.5truecm>
\putrule from -4 3 to 1.5 3
\setlinear  \plot  -1.33 2  4 6 /          \plot 2 6  5 2    /
            \plot  3 2  5.66 4 /         \plot 11 6  14 2 /
            \plot  10.33 4  13 6 /
            \plot  12 2  17.33 6 /        
 \setdashes  \plot   7 5  5.66 4  /       
             \plot   9 3  10.33 4 /      
\putrule from -6 3 to -4 3
 \setsolid
\put {\dots} at 8 4
\put {$E_r$} at 16.9 4.9
\put {$E_{r-1}$} at 13.4 4.5
\put {$E_{2}$} at 4 4.5
\put {$E_1$} at .9 4.5
\put {$E_i$} at -3 3.6
\endpicture}
\vskip .7truecm
\centerline{\smc Figure 1}
\vskip 1truecm

\noindent
{\bf 0.4.}  When $Y$ is a projective algebraic variety (of arbitrary dimension) with at worst Gorenstein canonical singularities, Batyrev also proposed in [B1] the following definition of stringy Hodge numbers $h^{p,q}_{st}(Y)$ {\sl if $E(Y)$ is a polynomial}.  Say $E(Y) = \sum_{p,q} b_{p,q} u^pv^q$ $\in \Bbb Z[u,v]$, then $h^{p,q}_{st} (Y) := (-1)^{p+q} b_{p,q}$.  He also conjectured that then all $h^{p,q}_{st}(Y)$ are nonnegative.  (When $Y$ is nonsingular $E(Y) = H(Y)$ and hence the stringy Hodge numbers are just the usual Hodge numbers.)

One can ask more generally for varieties with at worst log terminal singularities to define stringy Hodge numbers if $E(Y)$ is a polynomial (with rational powers).  Then however they would be rationally graded; now this phenomenon appears naturally in the stringy geometry of orbifolds, see Ruan [R].
\bigskip
\noindent
{\bf 0.5.}  Here we show that the stringy $E$--functions of log terminal surface singularities are polynomials in $uv$ (with rational powers) {\it with nonnegative coefficients}.
Hence it really makes sense to associate stringy Hodge numbers to a normal surface with at worst log terminal singularities via its stringy $E$--function.

This result uses a formula for the contribution of a chain of $\Bbb P^1$'s in $\cup_{i \in T} E_i$ to the stringy $E$--function in terms of the determinant of a certain non--symmetric deformation of the intersection matrix of the components in the chain,
see (5.5).  Such a determinant appeared already in [V1] and [V2].
\bigskip
\demo{0.6. Remark}  The stringy Euler numbers in this paper are different from the orbifold Euler numbers in the recent work of Langer [L]; as he already remarked stringy Euler numbers are topological invariants, which is not the case for orbifold Euler numbers. \enddemo
\bigskip
\noindent
{\bf 0.7.} In \S 1 we recall for the convenience of the reader some basic notions and the classification of log canonical surface singularities.  Theorem 0.3 is proved in \S 2.  Then in \S 3 we introduce our stringy invariants, working more generally with the Grothendieck ring of algebraic varieties instead of Hodge polynomials.  Generalizing the log terminal case [B1, Theorem 3.7], we verify in \S 4 Poincar\'e duality for the stringy $E$--function of a complete normal surface, and we describe the meaning of its constant term.  In \S5 we prove the nonnegativity of the coefficients of the determinant, expressing the contribution of a chain of $\Bbb P^1$'s to the stringy $E$--function; finally we use this in $\S 6$ to treat the log terminal singularities and to notice a remarkable fact concerning weighted homogeneous surface singularities.

\bigskip
\demo{Acknowledgement} We would like to thank B. Rodrigues for his useful remarks.
\enddemo

\bigskip
\bigskip
\heading
1. Basics
\endheading
\bigskip
\noindent
{\bf 1.1.} We denote by ${\Cal V}$ the Grothendieck ring of complex algebraic varieties (i.e. of reduced separated schemes of finite type over $\Bbb C$).  This is the free abelian group generated  by the symbols [V], where [V] is a variety, subject to the relations $[V] = [V^\prime]$ if $V \cong V^\prime$ and $[V] = [V \setminus W] + [W]$ if $W$ is closed in $V$.  Its ring structure is given by $[V]\cdot [W] := [V \times W]$.  We abbreviate $L := [\Bbb A^1]$.

For a variety $V$ we denote by $h^{p,q} (H^i_c(V, \Bbb C))$ the rank of the $(p,q)$--Hodge component in the mixed Hodge structure of the $i$th cohomology group with compact support of $V$, and we put $e^{p,q}(V) := \sum_{i \geq 0} (-1)^i h^{p,q} (H^i_c(V,\Bbb C))$.  The {\it Hodge polynomial} of $V$ is
$$H(V) = H(V;u,v) := \sum_{p,q} e^{p,q} (V) u^pv^q \in \Bbb Z[u,v] \, .$$
Precisely by the defining relations of ${\Cal V}$ there is a well defined ring morphism   \linebreak    $H : {\Cal V} \rightarrow \Bbb Z[u,v]$ determined by $[V] \mapsto H(V)$.

We denote by $\chi(V)$ the topological Euler characteristic of $V$, i.e. the alternating sum of the ranks of its Betti or the Rham cohomology groups.  Clearly $\chi(V) = H(V;1,1)$ and we also obtain a ring morphism $\chi : {\Cal V} \rightarrow \Bbb Z$ determined by $[V] \mapsto \chi (V)$.
\bigskip
\demo{1.2. Remark}  For a variety $Y$ as in (0.1) one can also consider the `finest' stringy invariant
$${\Cal E}(Y) := \sum_{I \subset T} [E^\circ_I] \prod_{i \in I} \frac{L-1}{L^{a_i}-1} ,$$
living in a ring which we will not specify here (see e.g. [V3]); it specializes via the map $H$ to $E(Y)$.
We will discuss this when $Y$ is a surface in (3.2).
\enddemo
\bigskip
\noindent
{\bf 1.3.} Let $P \in S$ be a normal complex surface germ and $\pi : X \rightarrow S$ a resolution of $P \in S$, i.e. a proper birational morphism from a smooth surface $X$.  This implies that $\pi^{-1} P = \cup_{i \in T} E_i$ where the $E_i$ are (distinct) irreducible curves.
We call $\pi$ a {\it log} (or {\it good}\,) resolution if moreover $\cup_{i \in T} E_i$ is a normal crossings divisor, meaning that the $E_i$ are smooth curves intersecting transversely.  Denoting by $K_{\centerdot}$ the canonical divisor, the expression
$$K_X = \pi^\ast K_S + \sum_{i \in T} (a_i - 1) E_i \tag 1$$
makes sense and $a_i \in \Bbb Q$ is called the log discrepancy of $E_i$.
\bigskip
\noindent
{\bf 1.4.}  For the sequel it is useful to recall here the construction of the $a_i$.  Although in arbitrary dimension it is not clear how to pullback Weil divisors (like $K_S$), there is a natural definition on normal surfaces by Mumford [M, p.17].  This amounts here in applying the adjunction formula to (1) for all $E_j$, yielding the linear system of equations
$$\sum_{i \in T} (a_i - 1) E_i \cdot E_j = 2p_a (E_j) - 2 - E^2_j, \quad j \in T, \tag 2$$
where $p_a(E_j)$ is the arithmetic genus of $E_j$.  Since the intersection matrix of the $E_i$ is negative definite [M], this linear system indeed has a unique solution over $\Bbb Q$.
\bigskip
\noindent
{\bf 1.5. Definition.}  The singularity $P \in S$ is called {\it canonical}, {\it log terminal} and {\it log canonical} if for some (or equivalently : any) log resolution of $P \in S$ we have that all $a_i, i \in T$, are $\geq 1, > 0$ and $\geq 0$, respectively.

The canonical (surface) singularities are precisely the Du Val or ADE singularities.
\bigskip
\noindent
{\bf 1.6.} We recall the classification of log terminal and log canonical surface singularities using the pictorial presentation of [A], to which we also refer for a proof and adequate references.  We use the dual graph of the {\it minimal} log resolution $\pi : X \rightarrow S$ of $P \in S$ in which the $E_i, i \in T$, are represented by dots and an intersection between them by a line connecting the corresponding dots.  Also here by an ellips we abbreviate a chain of arbitrary length $r \geq 1$ as in Figure 2,
where $E_1 \cong \cdots \cong E_r \cong \Bbb P^1$, and the number $n$ is the absolute value of the determinant of the intersection matrix of $E_1, \cdots , E_r$.
\vskip 1truecm
\centerline{
\beginpicture
\setcoordinatesystem units <.5truecm,.5truecm>
\ellipticalarc axes ratio 3:1 360 degrees from -4 0 center at -7 0
\put {$n$} at  -7  0
\put {$:=$} at -2 0
\putrule from 0 0 to 5 0
\putrule from 7 0 to 10 0
\put {\dots} at 6 0
\multiput {$\bullet$} at  0 0  2 0  4 0  8 0  10 0 /
\put {$E_1$} at 0 -.8
\put {$E_2$} at 2 -.8
\put {$E_3$} at 4 -.8
\put {$E_{r-1}$} at 8 -.8
\put {$E_r$} at 10 -.8
\endpicture}
\vskip .7truecm
\centerline{\smc Figure 2}
\vskip 1truecm

\noindent
In Figures 3 and 4 below all $E_i$ are rational, except in case (3).

(i) The log terminal surface singularities have possible dual graphs as in Figure 3.
Among these the canonical singularities are those where all $E_i$ have self--intersection number $-2$.  We also recall that case (1) consists of the Hirzebruch--Jung singularities $A_{n,q}$, where
$q$ can be taken as the absolute value of the determinant of the intersection matrix of either $E_1, \cdots , E_{r-1}$ or
$E_2,\cdots,E_r$ (see e.g. [BPV, III5]).

\vskip 1truecm
\centerline{
\beginpicture
\setcoordinatesystem units <.5truecm,.5truecm>
\ellipticalarc axes ratio 3:1 360 degrees from 5 10 center at 2 10
\put {$n$} at  2  10
\put {$(1)$} at  -5 10
\ellipticalarc axes ratio 3:1 360 degrees from 2 6 center at -1 6
\ellipticalarc axes ratio 3:1 360 degrees from 12 7.5 center at 9 7.5
\ellipticalarc axes ratio 3:1 360 degrees from 12 4.5 center at 9 4.5
\put {$n_1$} at  -1 6
\put {$n_2$} at  9 7.5
\put {$n_3$} at  9 4.5
\put {$(2)$} at  -5 6
\put {$\bullet$} at 4 6
\putrule from 2 6 to 4 6
\plot 6 4.5  4 6   6 7.5 /
\put{with $(n_1,n_2,n_3) = \left\{ \aligned & (2,2,n_3 \geq 2) \\ & (2,3,3) \\ & (2,3,4) \\ & (2,3,5) \endaligned \right.$} 
at 20 6
\endpicture}
\vskip .7truecm
\centerline{\smc Figure 3}
\vskip 1truecm
(ii) The log canonical surface singularities which are not log terminal, called {\it strictly log canonical} in the sequel, have   possible dual graphs as in Figure 4. Cases (3) and (4) are called simple elliptic and cusp singularities, respectively.
\vskip 1truecm
\centerline{
\beginpicture
\setcoordinatesystem units <.5truecm,.5truecm>
\put {$\bullet$} at  3  18
\put {$(3)$} at  -5 18
\put{elliptic curve} at 16 18
\multiput {$\bullet$} at 3 10  5.1 10.9  6 13  5.1 15.1  3 16  
                         0.9 15.1  0 13  0.9 10.9   /
\plot 3 10  5.1 10.9  6 13  5.1 15.1  3 16  
                         0.9 15.1  0 13  0.9 10.9  3 10  /
\put {$(4)$} at  -5 13
\put{a closed chain of length $r\geq 2$} at 19 13
\ellipticalarc axes ratio 3:1 360 degrees from 2 6 center at -1 6
\ellipticalarc axes ratio 3:1 360 degrees from 12 7.5 center at 9 7.5
\ellipticalarc axes ratio 3:1 360 degrees from 12 4.5 center at 9 4.5
\put {$n_1$} at  -1 6
\put {$n_2$} at  9 7.5
\put {$n_3$} at  9 4.5
\put {$(5)$} at  -5 6
\put {$\bullet$} at 4 6
\putrule from 2 6 to 4 6
\plot 6 4.5  4 6   6 7.5 /
\put{with $(n_1,n_2,n_3) = \left\{ \aligned & (2,3,6) \\ & (2,4,4) \\ & (3,3,3)  \endaligned \right.$}
at 19 6	
\ellipticalarc axes ratio 3:1 360 degrees from 6 1 center at 3 1
\multiput {$\bullet$} at  -2 2  -2 0  8 2  8 0 /
\plot -2 2  0 1 /           
\plot -2 0  0 1 /
\plot 8 2  6 1 /           
\plot 8 0  6 1 /
\put {$(6)$} at  -5 1	  
\put {where the four end curves }	 
at 18.5 1.5
\put {have self--intersection $-2$}	 
at 18.5 .5	 
\endpicture}
\vskip .7truecm
\centerline{\smc Figure 4}

\bigskip
\bigskip
\heading
2. Structure of the resolution graph
\endheading
\bigskip
\noindent
{\bf 2.1.}  Let $P \in S$ be a normal complex surface germ.  We exclude the easy case that $P \in S$ is canonical.

Let $\pi_0 : X_0 \rightarrow S$ be the minimal resolution of $P \in S$ and $\pi : X \rightarrow S$ its minimal log resolution.  So $\pi = h \circ \pi_0$ where $h : X \rightarrow X_0$ is the minimal embedded resolution of $\pi^{-1}_0 P \subset X_0$.  We will denote the irreducible components of $\pi^{-1}_0P$ and $\pi^{-1}P$ by $E_i, i \in T_0$ and $i \in T$, respectively (so $T_0 \subset T$).
For $i \in T$ we will denote by $a_i$ the log discrepancy of $E_i$.

\bigskip
\proclaim{2.2. Lemma {\rm[A, 3.1.2--3.1.3]}}
(i) $\forall i \in T_0 : a_i < 1$.

(ii) Let $T^\prime_0 \varsubsetneq T_0$.  Then the solutions $a^\prime_i, i \in T^\prime_0$, of the linear system (2'), which is given by replacing $T$ by $T^\prime_0$ in (2), satisfy $a_i < a^\prime_i$ for all $i \in T^\prime_0$. \endproclaim

\bigskip
\noindent
{\bf 2.3.}  Denote $\kappa_j = -E^2_j$ (in $X_0$) for $j \in T_0$; so $\kappa_j \in \Bbb N \setminus \{ 0 \}$, and if $E_j \cong
\Bbb P^1$, then $\kappa_j \geq 2$ because $\pi_0$ is the minimal resolution.  We will often use the equalities in (2) in the form
$$\kappa_j a_j = \sum^r_{i=1} (E_j \cdot E_i) (a_i - 1) + 2 - 2p_a (E_j) \tag 3$$
for $j \in T_0$, where $E_j$ intersects precisely $E_1,\cdots,E_r$.

We now start the discussion of the structure of the resolution graph $\cup_{i \in T_0} E_i \subset X_0$ with respect to the nullity and signs of the log discrepancies.  We begin with the following obvious remark.

\bigskip
\proclaim
{2.4. Remark}  If $p_a (E_j) \geq 1$ then either

(i) $a_j < 0$, or

(ii) $a_j = 0$, $p_a (E_j) = 1$ and $\pi^{-1}_0 P = E_j$.
\newline
More precisely there are three possibilities for $E_j$ in case (ii) : a (nonsingular) elliptic curve, a rational curve with one node or a rational curve with one (ordinary) cusp. \endproclaim

\medskip
\proclaim
{2.5. Lemma}  Suppose that $\cup_{i \in T_0} E_i$ is {\it not} a normal crossings divisor at the point $Q \in X_0$. Then either

(i) all $E_j \ni Q$ satisfy $a_j < 0$, or

(ii) we have that $a_i = 0$ for all $i \in T_0$, and more precisely $\cup_{i \in T_0} E_i$ is one of the following (see Figure 5)~:
\itemitem{(1)} a rational curve with one node,
\itemitem{(2)} a rational curve with one cusp,
\itemitem{(3)} two smooth rational curves intersecting in $Q$ with intersection multiplicity 2,
\itemitem{(4)} three smooth rational curves intersecting each other in $Q$ such that any two of them intersect transversely.
\endproclaim
\vskip .5truecm
\centerline{
\beginpicture
\setcoordinatesystem units <.5truecm,.5truecm> point at 1 0
\ellipticalarc axes ratio 2:3 80 degrees from 0 -1 center at 0 2
\circulararc 180 degrees from 0 1 center at 0 0
\ellipticalarc axes ratio 2:3 -80 degrees from 0 1 center at 0 -2
\put {$\bullet$} at 1.5 0
\put {$Q$} at  2.2 0
\put {$(1)$} at -2 0
\setcoordinatesystem units <.5truecm,.5truecm> point at -7 0
\ellipticalarc axes ratio 2:3 80 degrees from 0 0 center at 0 2
\ellipticalarc axes ratio 2:3 -80 degrees from 0 0 center at 0 -2
\put {$\bullet$} at 0 0
\put {$Q$} at  -.5 .5
\put {$(2)$} at -2 0
\setcoordinatesystem units <.5truecm,.5truecm> point at -14 0
\setquadratic   \plot 1 1.5   0 0   1 -1.5 /
                \plot -1 1.5   0 0   -1 -1.5 /
\put {$\bullet$} at 0 0
\put {$Q$} at  .6 0
\put {$(3)$} at  -2 0
\setcoordinatesystem units <.5truecm,.5truecm> point at -21 0
\putrule from -1.8 0 to 1.8 0
\setlinear
\plot -1 -1.5   1 1.5 /         \plot 1 -1.5   -1 1.5 /
\put {$\bullet$} at 0 0
\put {$Q$} at  .9 .5
\put {$(4)$} at -2.8 0
\endpicture}
\vskip 1truecm
\centerline{\smc Figure 5}
\vskip .7truecm

\demo
{Proof}  If $\pi^{-1}_0 P$ consists of just one component $E_j$, then by assumption this $E_j$ is not smooth and thus $p_a(E_j) \geq 1$.  Then by Remark 2.4 either $a_j < 0$ or we are dealing with case (1) or (2).

From now on we suppose that $|T_0| \geq 2$ and we denote $T^\prime_0 := \{ j \in T_0 \mid Q \in E_j \}$.  If $T^\prime_0 = \{ j \}$, then again $E_j$ must be singular, $p_a(E_j) \geq 1$ and $a_j < 0$. So we can assume further that also $|T^\prime_0| \geq
2$.  For $j \in T^\prime_0$ denote
$$c^\prime_j := 2p_a (E_j) - 2 + \sum_{\Sb k \ne j \\ k \in T^\prime_0 \endSb} E_k \cdot E_j.$$
Since $\cup_{k \in T^\prime_0} E_k$ is not a normal crossings divisor at $Q$, we have that each $c^\prime_j \geq 0$.  We consider the linear system of equations
$$\sum_{i \in T^\prime_0} a^\prime_i E_i \cdot E_j = c^\prime_j \quad , \quad j \in T^\prime_0 ,$$
and argue as in [A, 3.1.2].  Since the intersection matrix of the $E_i, i \in T^\prime_0$, is negative definite, all coefficients of its inverse matrix are strictly negative, and so all $a^\prime_i \leq 0$.
\bigskip
\noindent
{\sl First case :} $T^\prime_0  \varsubsetneq T_0$.  It should be clear that the $a^\prime_i, i \in T^\prime_0$, are also the solutions of the system (2') in Lemma 2.2(ii); hence this lemma yields that $a_i < 0$ for all $i \in T^\prime_0$.
\bigskip
\noindent
{\sl Second case :} $T^\prime_0 = T_0$.  So $a^\prime_i = a_i$ for $i \in T^\prime_0 = T_0$.  If at least one $c^\prime_j > 0$ then all $a_j < 0$.  On the other hand, if all $c^\prime_j = 0$, then for all $j \in T_0$ we have that $a_j = 0$, and by definition of $c^\prime_j$ also that $p_a(E_j) = 0$ and $\sum_{\Sb k \ne j \\ k \in T_0 \endSb} E_k \cdot E_j = 2$.  This is only possible when dealing with cases (3) or (4). \qquad $\square$  \enddemo

\bigskip
\proclaim{2.6. Proposition}  Suppose that $a_j \geq 0$ for all $j \in T_0$.  Then either

(i) $P \in S$ is log canonical, or

(ii) $P \in S$ has minimal resolution $\cup_{i \in T_0} E_i \subset X_0$ as in cases (2), (3) or (4) of Lemma 2.5, and thus minimal log resolution $\cup_{i \in T} E_i \subset X$ as in Figure 6, where all curves are rational and the numbers denote log discrepancies.
\endproclaim
\vskip .5truecm
\centerline{
\beginpicture
\setcoordinatesystem units <.5truecm,.5truecm> 
\putrule from 0 0 to 8 0
\putrule from 2 -1.5 to 2 1.5
\putrule from 4 -1.5 to 4 1.5
\putrule from 6 -1.5 to 6 1.5
\multiput {$0$} at  1.6 -1.6  3.6 -1.6  5.6 -1.6 / 
\put {$-1$} at 8.8 0.3
\endpicture}
\vskip .7truecm
\centerline{\smc Figure 6}
\pagebreak
\demo
{Proof}  If $\cup_{i \in T_0} E_i$ is a normal crossings divisor, then $X_0 = X$ and $P \in S$ is log canonical.  Otherwise we have one of the four cases of Lemma 2.5(ii), where case (1) is also log canonical. \quad $\square$ \enddemo
\bigskip
\noindent
The following lemma will be the crucial ingredient for our structure theorem.

\smallskip
\proclaim
{2.7. Lemma}  Let a fixed $E_j, j \in T_0,$ intersect (in $X_0$) another component $E_1$ such that their log discrepancies satisfy
$a_j \geq 0$  and $a_1 < a_j$.
Then

(i) $E_j \cong \Bbb P^1$ and $\cup_{i \in T_0} E_i$ is a normal crossings divisor at {\it each} point of $E_j$; also

(ii) $E_j$ intersects {\it at most twice} $\cup_{\ell \neq j} E_\ell$, and if it intersects besides $E_1$ also $E_i$, then $E_i \cong \Bbb P^1$ and $a_j < a_i$.
\endproclaim
\vskip .5truecm
\centerline{
\beginpicture
\setcoordinatesystem units <.5truecm,.5truecm>
\putrule from -4 3 to 1.5 3
\setlinear  \plot  -1.33 2  4 6 /   
\setdashes     \putrule from -6 3 to -4 3
 \setdashes <3pt>  
 \plot 2 6  5 2    /
\put {$E_{i}$} at 4.9 3.7
\put {$E_j$} at .9 4.5
\put {$E_1$} at -3 3.6
\endpicture}
\vskip .7truecm
\centerline{\smc Figure 7}
\vskip 1truecm

\demo
{Proof} (i) This is immediate by Remark 2.4 and Lemma 2.5, respectively.

(ii) Suppose that $E_j$ intersects besides $E_1$ at least twice   other components, say $E_2$ and $E_3$.  The equation (3) for $E_j$ yields
$$\kappa_ja_j \leq \sum^3_{i=1} (a_i - 1) + 2 = a_1 + a_2 + a_3 - 1 \, ,$$
since possible other components $E_\ell$ intersecting $E_j$ have $a_\ell < 1$.  Using furthermore that $a_1 < a_j$ we obtain
$$0 \leq (\kappa_j - 1)a_j < a_2 + a_3 - 1 . \tag 4$$
Since also $a_2 < 1$ and $a_3 < 1$ this can only occur if both $a_2 > 0$ and $a_3 > 0$.  Then by Remark 2.4 we have $E_2 \cong E_3 \cong \Bbb P^1$, and thus $\kappa_2, \kappa_3 \geq 2$.  Applying now equation (3) to $E_i$ for $i = 2,3$ yields $\kappa_i a_i \leq a_j + 1$ and thus
$$a_i \leq \frac{a_j +1}{\kappa_i} \leq \frac{a_j+1}{2} \, .$$
Plugging in these inequalities in (4) we obtain
$$(\kappa_j - 1)a_j < a_j \, ,$$
which is impossible since $\kappa_j \geq 2$ and $a_j \geq 0$.  So indeed $E_j$ intersects besides $E_1$ at most once another component.

Suppose finally that also $E_i$ intersects $E_j$.  Then $\kappa_j a_j = a_1 + a_i$, implying analogously $a_j \leq (\kappa_j - 1)a_j < a_i$.  So certainly $a_i > 0$ and hence $E_i \cong \Bbb P^1$. \qquad $\square$ \enddemo

\bigskip
\proclaim
{2.8. Theorem}  Suppose that at least one $E_j, j \in T_0$, has log discrepancy $a_j < 0$.  Then the minimal resolution graph  $\pi^{-1}_P = \cup_{i \in T_0} E_i (\subset X_0)$ of $P \in S$ consists of the connected part ${\Cal N}_0 := \cup_{\Sb i \in T_0 \\ a_i < 0 \endSb} E_i$, to which a finite number of chains
are attached as in Figure 1. Here $E_i \subset {\Cal N}_0, E_\ell \cong \Bbb P^1$ for $1 \leq \ell \leq r$, $a_1 \geq 0$ and $(a_i <) a_1 < a_2 < \cdots < a_r$, and $\cup_{i \in T_0} E_i$ is a normal crossings divisor at each point of $\cup^r_{\ell=1} E_\ell$.
%
\endproclaim

\medskip
\demo
{Proof}  If $a_i < 0$ for all $i \in T_0$, then of course ${\Cal N} = \cup_{i \in T_0} E_i$ is connected and there is no other claim.

On the other hand suppose that there is an $E_\ell, \ell \in T_0$, with $a_\ell \geq 0$; then there is also a component $E_1$ with $a_1 \geq 0$ that intersects a component of ${\Cal N}$, say $E_i$.  Now by Lemma 2.7 either $E_1$ intersects only $E_i$, or
it intersects exactly one other component $E_2$ with $a_1 < a_2$ (see Figure 8).
\vskip 1truecm
\centerline{
\beginpicture
\setcoordinatesystem units <.5truecm,.5truecm>
\putrule from -4 3 to 1.5 3
\setlinear  \plot  -1.33 2  4 6 /   
 \setdashes     \putrule from -6 3 to -4 3
 \setdashes <3pt>  
 \plot 2 6  5 2    /
\put {$E_{2}$} at 4.9 3.7
\put {$E_1$} at .9 4.5
\put {$E_i$} at -3 3.6
\endpicture}
\vskip .7truecm
\centerline{\smc Figure 8}
\vskip 1truecm
\noindent
In the last case analogously either $E_2$ intersects only $E_1$, or exactly one other component $E_3$ with $a_2 < a_3$.  Ultimately we must obtain a chain as described above, where the rationality and normal crossings property of its components are implied by Lemma 2.7.

Now all other $E_k$ with $a_k \geq 0$ form such (necessarily disjoint) chains, implying that ${\Cal N_0}$ is connected. \qquad $\square$ \enddemo
\bigskip
\demo{2.9. Addendum (easy exercise)}  Suppose that $a_1 = 0$ in Theorem 2.8.

\itemitem{(i)} If $r = 1$, then $a_i = -1$.
\itemitem{(ii)} If $r = 2$, then $a_i = - \frac{1}{\kappa_2} \geq -\frac 12$.
\itemitem{(iii)} If $r > 2$, then $a_i > -\frac{1}{\kappa_2} \geq -\frac 12$.
\enddemo
\vskip 1truecm
\centerline{
\beginpicture
\setcoordinatesystem units <.5truecm,.5truecm>
\putrule from -4 3 to 1.5 3
\setlinear  \plot  -1.33 2  4 6 /          \plot 2 6  5 2    /
            \plot  3 2  5.66 4 /         \plot 11 6  14 2 /
            \plot  10.33 4  13 6 /
            \plot  12 2  17.33 6 /        
 \setdashes  \plot   7 5  5.66 4  /       
             \plot   9 3  10.33 4 /      
\putrule from -6 3 to -4 3
 \setsolid
\put {\dots} at 8 4
\put {$E_r$} at 16.9 4.9
\put {$E_{r-1}$} at 13.4 4.5
\put {$E_{2}$} at 4 4.5
\put {$E_1$} at .9 4.5
\put {$E_i$} at -3 3.6
\endpicture}
\vskip .7truecm
\centerline{\smc Figure 1}
\vskip 1truecm
\noindent
{\bf 2.10.} We now derive from Theorem 2.8 precisely the same statement for the minimal {\sl log} resolution of any singularity which is not log canonical.  For completeness we recall the data.

\bigskip
\proclaim
{Theorem}  Let $P \in S$ be a normal surface singularity germ which is not log canonical.  Let $\pi : X \rightarrow S$ be the minimal log resolution of $P \in S$; denote the irreducible components of $\pi^{-1}P$ by $E_i, i \in T$, and their log discrepancies by $a_i$.  Then $\pi^{-1} P = \cup_{i \in T} E_i$ consists of the connected part ${\Cal N} = \cup_{\Sb i \in T \\ a_i < 0 \endSb} E_i$,  to which a finite number of chains
are attached as in Figure 1. Here $E_i \subset {\Cal N}, E_\ell \cong \Bbb P^1$ for  $1 \leq \ell \leq r$, $a_1 \geq 0$ and $(a_i <) a_1 < a_2 < \cdots < a_r <1$.

In particular either all $a_\ell > 0$ for $1 \leq \ell \leq r$, or $a_1 = 0$ and $a_\ell > 0$ for $2 \leq \ell \leq r$. \endproclaim

\medskip
\demo
{Proof}  By Proposition 2.6 either $\cup_{i \in T} E_i \subset X$ is as in Figure 6, which is as stated, or at least one $E_i, i \in T_0$, has log discrepancy $a_i < 0$ and then we can apply Theorem 2.8.  The description of $\cup_{i \in T_0} E_i$ there implies the following.  Let $Q \in \cup_{i \in T_0} E_i$ be the centre of one of the first blowing--ups of $h : X \rightarrow X_0$ (if $X \ne X_0$).  Then {\it all} $E_\ell$ that contain $Q$ belong to ${\Cal N}_0$, i.e. have $a_\ell < 0$.

We only have to show that the exceptional curve $E$ of this blowing--up still has log discrepancy
$a < 0$, because then $\cup_{i \in T}E_i$ will be as stated by repeated application of this argument.

Now indeed the equality (3) applied to $E$, considered in the blown--up of $X_0$ in $Q$, yields
$$a = \sum_{\Sb \ell \in T_0 \\ E_\ell \ni Q \endSb} \mu_\ell (a_\ell - 1) + 2 \, ,$$
where $\mu_\ell$ is the multiplicity of $E_\ell \subset X_0$ in the point $Q$.  Since $\cup_{\Sb \ell \in T_0 \\ E_\ell \ni Q \endSb} E_\ell$ is not a normal crossings divisor in $Q$, we have that
$\sum_{\Sb \ell \in T_0 \\ E_\ell \ni Q \endSb} \mu_\ell \geq 2$ and consequently $a < 0$.

Remark that thus all $a_\ell, \ell \in T,$ satisfy $a_\ell <1$, which was not a priori obvious for the minimal {\sl log} resolution of
$P \in S$.
\qed
\enddemo

\bigskip
\noindent
{\bf 2.11.} We will not use it in this paper, but Theorem 2.10 can easily be extended to the situation $P \in (S,B)$, where $B$ is a reduced divisor on $S$ such that each of its irreducible components contains $P$. For any resolution $\tilde \pi : \tilde
X \rightarrow S$ of $P \in S$, denoting the irreducible components of ${\tilde \pi}^{-1}B$ by $E_i, i\in \tilde T,$ the log discrepancies $a_i, i \in \tilde T,$ are now defined by $K_{\tilde X}={\tilde \pi}^*(K_S+B)+\sum_{i\in \tilde T} (a_i -1)E_i$. Note that in particular $a_i=0$ for each component $E_i$ of the strict transform of $B$.

Working through the arguments in this section, one can verify that Theorem 2.10 remains true, where now $\pi$ is the minimal log resolution of (the germ at $P$ of) the pair $(S,B)$ and $T$ runs over the irreducible components of $\pi^{-1}B$. In particular the irreducible components of the strict transform of $B$ give \lq chains with $r=1$ and $a_1=0$\rq.

\bigskip
\bigskip
\heading
3. Definition of stringy invariants
\endheading
\bigskip
\noindent
{\bf 3.1.}  Here we let $S$ be a normal algebraic surface.  We denote again by $\pi : X \rightarrow S$ the minimal log resolution of $S$ and by $E_i, i \in T$, the irreducible components of the exceptional divisor of $\pi$ with log discrepancies $a_i, i \in T$.  For $I \subset T$ we put $E_I := \cap_{i \in I} E_i$ and $E^\circ_I := E_I \setminus \cup_{\ell \not\in I} E_\ell$.  (So $E_\emptyset = X$.)
\bigskip
\noindent
{\bf 3.2.} It is clear that when {\it all} log discrepancies $a_i \ne 0$ one can at least define stringy invariants $e(S)$ and $E(S)$ by the same formulas as in (0.1).  One can moreover consider the `finest' such invariant ${\Cal E}(S)$ on the level of the Grothendieck ring ${\Cal V}$ of algebraic varieties.
Denote by $t$ the least common denominator of the $a_i, i \in T$, and by $L^{1/t}$ the class of $Y$ in the quotient of the polynomial ring ${\Cal V}[Y]$ by $(Y^t - L)$.  Then
$${\Cal E}(S) = \sum_{I \subset T} [E^\circ_I] \prod_{i \in I} \frac{L-1}{L^{a_i} - 1},$$
living in the localization of ${\Cal V}[L^{1/t}, L^{-1/t}]$ with respect to the elements $L^{b/t} - 1,$ \linebreak  $b \in \Bbb Z \setminus \{ 0 \}$.  One can easily generalize $H : {\Cal V} \rightarrow \Bbb Z[u,v]$ to a ring morphism $H$ from this ring to the `rational functions in $u,v$ with rational powers' such that $H({\Cal E}(S)) = E(S)$.

We now extend in a natural way these notions to surfaces $S$ for which certain log discrepancies are allowed to be zero; more precisely we allow $a_i = 0$ if $E_i \cong \Bbb P^1$ and $E_i$ intersects exactly once or twice other components $E_\ell$ (with $a_\ell \neq 0$).  By Theorem 2.10 this means that we consider in fact {\it all} normal surfaces without strictly log canonical singularities !
\bigskip
\noindent
{\bf 3.3.}  We first motivate our definition.  Let $E_i \cong \Bbb P^1$ intersect exactly $E_1$ and $E_2$, and denote $E^2_i = - \kappa_i$.  If all log discrepancies are nonzero, then the contribution of $E_i$ to (the generalized) $e(S)$ is
$$\frac{0}{a_i} + \frac{1}{a_ia_1} + \frac{1}{a_ia_2} = \frac{a_1+a_2}{a_ia_1a_2} = \frac{\kappa_i}{a_1a_2}\, .$$
This expression is also meaningful if $a_i = 0$, and can thus be considered as generalizing the `classical' contribution of $E_i$ if $a_i \ne 0$.  (Analogously when $E_i$ intersects only $E_1$ its contribution is $\frac{\kappa_i}{a_1}$.)

In the same spirit (when all log discrepancies are nonzero) the contribution of $E_i$ to ${\Cal E}(S)$ is
$$ \split & (L-1) \frac{L-1}{L^{a_i}-1} + \frac{(L-1)^2}{(L^{a_i}-1)(L^{a_1}-1)} + \frac{(L-1)^2}{(L^{a_i}-1)(L^{a_2}-1)} \\ & = \frac{(L-1)^2(L^{\kappa_ia_i}-1)}{(L^{a_i}-1)(L^{a_1}-1)(L^{a_2}-1)}  = \frac{(L-1)^2(L^{(\kappa_i-1)a_i}+ \cdots + L^{a_i}+1)}{(L^{a_1}-1)(L^{a_2}-1)}\, . \endsplit$$
Also this expression makes sense if $a_i = 0$ and becomes then
$$\frac{\kappa_i(L-1)^2}{(L^{a_1}-1)(L^{a_2}-1)}\, .$$

\bigskip
\proclaim{3.4. Definition} \rm  Let $S$ be a normal algebraic surface without strictly log canonical singularities for which we use the notation of (3.1). We put furthermore $Z := \{ i \in T \mid a_i = 0 \}$.

(i) The {\it stringy Euler number} of $S$ is
$$e(S) :=  \sum_{I \subset T \setminus Z} \chi (E^\circ_I) \prod_{i \in I} \frac{1}{a_i} + \sum_{i \in Z} \frac{\kappa_i}{a_{i_1} a_{i_2}} \, ,$$
where for $i \in Z$ we denote $-\kappa_i = E^2_i$, and $E_i$ intersects either $E_{i_1}$ and $E_{i_2}$ or only $E_{i_1}$ (and then we put $a_{i_2} := 1$).

(ii) With the same notation the {\it stringy ${\Cal E}$--invariant} and {\it stringy $E$--function} of $S$ are
$${\Cal E}(S) := \sum_{I \subset T \setminus Z} [E^\circ_I] \prod_{i \in I} \frac{L-1}{L^{a_i}-1} + \sum_{i \in Z} \frac{\kappa_i(L-1)^2}{(L^{a_{i_1}}-1)(L^{a_{i_2}}-1)}  \qquad\text{and}$$
$$E(S) := \sum_{I \subset T \setminus Z} H(E^\circ_I) \prod_{i \in I} \frac{uv-1}{(uv)^{a_i}-1} + \sum_{i \in Z} \frac{\kappa_i(uv - 1)^2}{((uv)^{a_{i_1}}-1)((uv)^{a_{i_2}}-1)}.$$
So as before $H({\Cal E}(S))) = E(S)$ and $\lim_{u,v \rightarrow 1} E(S) = e(S)$. \endproclaim
\bigskip
\noindent
{\bf 3.5.} The following supports the naturality of our definition : the stringy invariants of Definition 3.4 can in fact be given by the same expressions in terms of {\it any} log resolution with allowed zero discrepancies as above.

\medskip
\proclaim
{Proposition}  Let $S$ be a normal algebraic surface without strictly log canonical singularities.  Let $\pi^\prime : X^\prime \rightarrow S$ be a log resolution of $S$ and $E_i, i \in T^\prime (\supset T)$, the irreducible components of the exceptional divisor of $\pi^\prime$ with log discrepancies $a_i, i \in T^\prime$.  Denote $Z^\prime := \{ i \in T^\prime \mid a_i = 0 \}$.

We suppose that if $a_i = 0$ for $i \in T^\prime$, then $E_i$ intersects at most twice other components.  Then the invariants $e(S), E(S)$ and ${\Cal E}(S)$ of (3.4) are given by the same formula, using $T^\prime$ and $Z^\prime$ instead of $T$ and $Z$. \endproclaim

\medskip
\demo
{Proof}  We may restrict ourselves to the case that $h : X^\prime \rightarrow X$ is the blowing--up of $X$ in a point $Q$.  Denote the exceptional curve of $h$ by $E$ with log discrepancy $a$.  We must only check what happens with zero log discrepancies;
there are essentially five cases to consider.
\medskip
\itemitem{(1)} $\{ Q \} = E_1 \cap E_2$, $a_1 < 0 < a_2$ and $a = a_1 + a_2 = 0$,
\itemitem{(2)} $Q \in E^\circ_1$ and $a = a_1 + 1 = 0$,
\itemitem{(3)} $E_1 \cong \Bbb P^1$ has $a_1 = 0$, $E_1$ intersects $E_2$ and $E_3$ and $Q = E_1 \cap E_2$ (so $a_2 + a_3 = 0$ and $a  = a_2$),
\itemitem{(4)} $E_1 \cong \Bbb P^1$ has $a_1 = 0$, $E_1$ intersects only $E_2$ and $Q = E_1 \cap E_2$ (so $a = a_2 = -1$),
\itemitem{(5)} $E_1 \cong \Bbb P^1$ has $a_1 = 0$, $E_1$ intersects only $E_2$ and $Q \in E^\circ_1$ (so $a_2 = -1$ and $a = 1$).
\vskip 1truecm
\centerline{
\beginpicture
\setcoordinatesystem units <.5truecm,.5truecm>
\putrectangle corners at -2.5 2.5 and 6.5 -3.5
\putrule from -1.5 0 to 5.5 0
\putrule from 0 1.5 to 0 -2
\putrule from 4 1.5 to 4 -2
\setdots <3pt>
\putrule from 0 -2 to 0 -3
\putrule from 4 -2 to 4 -3
\setsolid
\put {$\bullet$} at 0 0
\put {$E_1$} at 2 .5
\put {$E_2$} at -.7 -1.5
\put {$E_3$} at 4.7 -1.5
\put {$Q$} at -.6 .6
\put {$X$} at -3.5 -.5
\put {$\longleftarrow$} at 9 -.5
\put {$h$} at 9 0
\setcoordinatesystem units <.5truecm,.5truecm> point at -15.5 .75
\putrectangle corners at -4 4 and 6.5 -3.5 
\putrule from 4 2 to 4 -2
\plot 5 0  0 3 /
\plot 2 3  -3 0 /
\plot -2 2  -1 -2 /
\setdots <3pt>
\putrule from 4 -2 to 4 -3
\plot -1 -2 -.75 -3 /
\setsolid
\put {$E_1$} at 2.3 .9
\put {$E_2$} at -1.8 -1.5
\put {$E_3$} at 4.7 -1.5
\put {$E$} at -.4 .9
\put {$X'$} at 7.5 .25
\endpicture}
\vskip .7truecm
\centerline{\smc Figure 9}
\vskip 1truecm
\noindent
We compute for example case (3) for ${\Cal E}(S)$; see Figure 9. We must verify whether the contributions to ${\Cal E}(S)$ of $E_1 \subset X$ and $E_1 \cup E \subset X^\prime$ are the same.  Let $-\kappa_1$ denote the self--intersection number of $E_1$ on
 $X$; then these contributions are
$$\frac{\kappa_1(L-1)^2}{(L^{a_2}-1)(L^{a_3}-1)}\quad \text { and } \quad \frac{(\kappa_1 + 1)(L-1)^2}{(L^a-1)(L^{a_3}-1)} + \frac{(L-1)^2}{L^a-1} + \frac{(L-1)^2}{(L^a-1)(L^{a_2}-1)}\, ,$$
respectively.  Since $a = a_2$ we have to show that
$$\frac{(L-1)^2}{L^{a_2}-1} \big(\frac{1}{L^{a_3}-1} + 1 + \frac{1}{L^{a_2}-1}\big) = \frac{(L-1)^2}{L^{a_2}-1} \cdot \frac{L^{a_2+a_3}-1}{(L^{a_3}-1)(L^{a_2}-1)}$$
is zero.  Indeed $a_2 + a_3 = 0$ here.

We leave the other cases as an exercise.  See [ACLM, Proposition 2.22] for similar computations. \qed
\enddemo

\bigskip
\demo{3.6. Remark}
Considering Theorem 2.10, the stringy invariants are in fact defined precisely for the normal algebraic surfaces $S$ for which the log discrepancies of the exceptional curves on the relative log minimal model of the pair $(S,0)$ are nonzero. See e.g. [KM][KMM] for this notion.
\enddemo

\bigskip
\demo{3.7. Remark}  Letting $t$ denote the least common denominator of the $a_i, i \in T$, it is clear that $E(S) \in \Bbb Q(u^{1/t}, v^{1/t}) \cap \Bbb Z[[u^{1/t},v^{1/t}]]$.  Also one can easily verify that the degree in $u$ or $v$ of $E(S)$ (as rational function) is at most 2.
\enddemo
\bigskip
\demo{3.8. Stringy invariants of germs}  When $P \in S$ is a normal surface germ, using the same notation as in (3.4) we define
$$e_P(S) := \sum_{\emptyset \ne I \subset T \setminus Z} \chi(E^\circ_I) \prod_{i \in I} \frac{1}{a_i} + \sum_{i \in Z} \frac{\kappa_i}{a_{i_1} a_{i_2}} \, ,$$
and analogously $E_P(S)$ and ${\Cal E}_P(S)$.  Remark 3.7 also applies to $E_P(S)$. \enddemo
\bigskip
\bigskip
\heading
4. Poincar\'e duality
\endheading
\bigskip
\noindent
{\bf 4.1.} Still using the notation of (3.1) we will generalize the Poincar\'e duality result of Batyrev [B1, Theorem 3.7] in the log terminal case to all complete normal surfaces $S$ for which our stringy invariants were defined, i.e. without strictly log canonical singularities.

When all log discrepancies $a_i \ne 0$ one sees easily that an alternative expression for ${\Cal E}(S)$ is
$${\Cal E}(S) = \sum_{I \subset T} [E_I] \prod_{i \in I} \big(\frac{L-1}{L^{a_i}-1} - 1\big).$$
For $E(S)$ this is the essential ingredient of Batyrev's proof.  (A similar expression and result appeared first in [DM].)  An analogous expression is true for our more general case.
\bigskip
\noindent
\proclaim
{4.2. Lemma}  Let $S$ be a normal algebraic surface without strictly log canonical singularities for which we use the notation of (3.4).  Then
$${\Cal E}(S) = \sum_{I \subset T \setminus Z} [E_I] \prod_{i \in I} \big(\frac{L-1}{L^{a_i}-1}-1\big) + \sum_{i \in Z}  \frac{\kappa_i (L-1)^2}{(L^{a_{i_1}}-1)(L^{a_{i_2}} - 1)}\, .$$    \endproclaim
\medskip
\demo
{Proof}  We must verify that
$$\sum_{I \subset T \setminus Z} [E_I] \prod_{i \in I} \big(\frac{L-1}{L^{a_i}-1}-1\big) = \sum_{I \subset T \setminus Z} [E^\circ_I] \prod_{i \in I} \frac{L-1}{L^{a_i}-1} \, .$$
Remembering that $[E_\emptyset] = [X] = \sum_{I \subset T} [E^\circ_I]$ one can compute that this amounts to the equality
$$\sum\Sb i \in Z \\ a_{i_2}\neq 1\endSb \left( [E^\circ_i] + \frac{L-1}{L^{a_{i_1}}-1} + \frac{L-1}{L^{a_{i_2}}-1} \right)
+\sum\Sb i \in Z \\ a_{i_2}= 1\endSb \left( [E^\circ_i] + \frac{L-1}{L^{a_{i_1}}-1} \right)
= 0.$$
This is true because precisely for $i \in Z$ we have in the first sum $[E^\circ_i] = L-1$ and $a_{i_1} + a_{i_2}$ (= $\kappa_i a_i$) = $0$, and in the second sum $[E^\circ_i] = L$ and $a_{i_1} = -1$ ! \qed  \enddemo

\bigskip
\proclaim
{4.3.  Poincar\'e duality theorem}  Let $S$ be a complete normal (algebraic) surface without strictly log canonical singularities.  Then
$$E(S;u,v) = (uv)^2 E(S;u^{-1},v^{-1}).$$ \endproclaim
\medskip
\demo
{Proof}  With the previous notation we have that
$$E(S;u,v) = \sum_{I \subset T \setminus Z} H(E_I) \prod_{i \in I} \big(\frac{uv-1}{(uv)^{a_i}-1} - 1\big) + \sum_{i \in Z} \frac{\kappa_i(uv - 1)^2}{((uv)^{a_{i_1}}-1)((uv)^{a_{i_2}} - 1)} \, .$$
The required functional equation is true for each term in the first sum (see [B1]; one should just remark that $E_\emptyset = X$ is projective because it is nonsingular, complete and of dimension 2).  And since for $i \in Z$ we have that $a_{i_1} + a_{i_2} = 0$, this functional equation is also true for each term in the last sum. \qed \enddemo
\bigskip
\noindent
{\bf 4.4.} When the projective surface $S$ has at worst log terminal singularities, one sees immediately that $E(S;0,0) = 1$ [B1, Theorem 3.7(ii)].  For more general $S$ this invariant measures as follows the complexity of the dual graphs of the non log terminal singularities of $S$.

\bigskip
\proclaim
{4.5. Proposition}   Let $S$ be a complete normal (algebraic) surface without strictly log canonical singularities.  Then
$$E(S;0,0) = 1 - \sum_{P_j} \chi (\Gamma_j) \, ,$$
where $P_j$ runs over the non log terminal singularities of $S$, and $\Gamma_j$ is the dual graph of the minimal (or any other) log resolution of $P_j \in S$. \endproclaim

\medskip
\demo
{Remark}  Recall that the first Betti number $\beta_1(\Gamma_j)$ of $\Gamma_j$ is $1 - \chi(\Gamma_j)$, and can be considered as the `number of circles' in $\Gamma_j$. \enddemo

\medskip
\demo
{Proof}  We use the expression for $E(S)$ coming from Lemma 4.2, where for $i \in Z$ we may assume that $a_{i_2} = -a_{i_1} > 0$ :
$$E(S) = H(X) + \sum_{\emptyset \ne I \subset T \setminus Z} H(E_I) \prod_{i \in I} \big(\frac{uv-1}{(uv)^{a_i}-1} - 1\big) - \sum_{i \in Z} \frac{\kappa_i (uv-1)^2(uv)^{a_{i_2}}}{((uv)^{a_{i_2}}-1)^2}.$$
Now for $i \not\in Z$, substituting $u = v = 0$ in $\frac{uv-1}{(uv)^{a_i}-1} - 1$ yields $0$ and $-1$ if $a_i > 0$ and $a_i < 0$, respectively.  So
$$\split E(S;0,0) & = H(X;0,0) +  \sum\limits_{\Sb \emptyset \ne I \subset T \\ \forall i \in I : a_i < 0 \endSb} H(E_I;0,0)(-1)^{|I|} \\ & = 1 - \dsize \sum\limits_{\Sb \emptyset \ne I \subset T \\ \forall i \in I : a_i < 0 \endSb} (-1)^{|I|+1}. \endsplit$$
The last sum can be interpreted as follows.  Looking at the terms associated to one fixed (necessarily non log terminal) singularity $P_j$, we see that they add to
\medskip
\centerline{(number of irreducible components $E_i$ with $a_i < 0$ in the minimal log resolution of $P_j$)}
\smallskip
\centerline{$-$ (number of intersections of such curves).}
\medskip
Using the structure of the resolution graph of $P$, as given by Theorem  2.10, we see that this difference is precisely the Euler characteristic $\chi(\Gamma_j)$.  Indeed, since the components $E_i$ with $a_i \geq 0$ appear in the disjoint `attached chains', they do not contribute to this Euler  characteristic. \qed \enddemo
\bigskip
\bigskip
\heading
5. Nonnegativity  of `stringy Hodge numbers'
\endheading
\bigskip
\noindent
{\bf 5.1.}  For projective varieties $V$ (of arbitrary dimension) with at worst canonical Gorenstein singularities Batyrev [B1, 3.10] conjectured the following.  If $E(V)$ is a polynomial, say $E(V) = \sum_{i,j} b_{ij} u^i v^j \in \Bbb Z[u,v]$, then all $(-1)^{i+j} b_{ij}$ should be nonnegative.  Hence it would make sense to consider them as {\it stringy Hodge numbers} of $V$ [B1, 3.8].

One can ask more generally an analogous question for arbitrary log terminal singularities.  Then however we talk about polynomials in $u,v$ with rational powers and hence about rationally graded stringy Hodge numbers.  We remark here that in the orbifold
philosophy of Ruan [R] such a rational grading of Hodge numbers appears in a natural way.

It is not totally obvious which nonnegativity to expect.  One can for instance consider $E(V)$ as a (usual) rational function in $u^{1/r}$ and $v^{1/r}$, where $r$ is the index of $V$.  At any rate when $E(V) = \sum_{i,j} b_{ij} u^i v^j$ is a polynomial over $\Bbb Z$ with rational powers, all $b_{ii}$ should be nonnegative in a decent theory of such stringy Hodge numbers.
\bigskip
\noindent
{\bf 5.2.}  In dimension 2 the canonical (Gorenstein) singularities are precisely those with all log discrepancies equal to 1 (in the minimal log resolution), hence their stringy $E$--function is by definition a polynomial with the required nonnegativity
of its coefficients.  More generally for log terminal surface singularities $P \in S$ one can easily compute that their stringy Euler number $e_P(S) \in \Bbb Z$ (in fact $\in \Bbb N$), as was also remarked by Batyrev [B1, 5.4].  Here we will express their stringy $E$--functions $E_P(S)$ in terms of the determinants of certain non--symmetric matrices which appeared already in [V1,V2], and derive
that $E_P(S)$ is a polynomial (with rational powers) in $uv$ with nonnegative coefficients.  Consequently normal surfaces with at worst log terminal singularities have well defined (rationally graded) stringy Hodge numbers in the sense of [B1].
\bigskip
\noindent
{\bf 5.3.} Let first $S$ be any normal surface without strictly log canonical singularities for which we still use the notation of 3.1.  We make a remark about obstructions for $E(S)$ to be a polynomial (with rational powers).  Looking at the defining formula for $E(S)$ each component $E_i$ with $a_i \ne 0$ is a possible obstruction because of the occurring expression $\frac{uv-1}{(uv)^{a_i}-1}$.  (However in the special case that $a_i \in \frac{1}{\Bbb Z}$ this expression is a polynomial.)

Now it is already implicit in [V2, 5.6] that chains of $\Bbb P^1$'s in fact do not contribute to the possible denominator of $E(S)$ (without special information on their log discrepancies), and more precisely we have a closed formula for the contribution
of such chains to $E(S)$.  We will state this more generally for ${\Cal E}(S)$.

\bigskip
\proclaim
{5.4. Proposition}  Let $E_i, 0 \leq i \leq r+1,$ intersect as in Figure 10, with $E_i \cong \Bbb P^1$ for $1 \leq i \leq r$ and $a_0 a_{r+1} \ne 0$.  Then the contribution of $\cup^r_{i=1} E_i$ to ${\Cal E(S)}$, i.e. the contribution of all $I \subset T$ with $I \cap \{ 1, \cdots , r \} \ne \emptyset$ is
$$\frac{(L-1)^2D_r}{(L^{a_0}-1)(L^{a_{r+1}}-1)}\, , \tag 5$$
where $D_r$ is the determinant of a non--symmetric `L--deformation' of the intersection matrix of the $E_i, 1 \leq 1 \leq r$, defined below.

For the same intersection configuration without $E_0$, and without $E_0$ and $E_{r+1}$, the contribution of $\cup^r_{i=1} E_i$ to ${\Cal E}(S)$ is $\frac{(L-1)D_r}{L^{a_{r+1}} - 1}$ and $D_r$, respectively.  (The last case means that $\{ 1, \cdots , r \}= T$.)
\endproclaim
\vskip .5truecm
\centerline{
\beginpicture
\setcoordinatesystem units <.5truecm,.5truecm>
\putrule from 0.5 1.5 to 0.5 6
\putrule from 16 1.5 to 16 6
\setlinear  \plot  -.5 3  3.5 6 /          \plot 2 6  5 2    /
            \plot  3 2  5.66 4 /         \plot 11 6  14 2 /
            \plot  10.33 4  13 6 /
            \plot  12 2  17.33 6 /
 \setdashes  \plot   7 5  5.66 4  /
             \plot   9 3  10.33 4 /
\putrule from 0.5 0 to 0.5 2
\putrule from 16 0 to 16 2
 \setsolid
\put {\dots} at 8 4
\put {$E_0$} at 16.8 .5
\put {$E_1$} at 14.5 4.7
\put {$E_2$} at 12 3.6
\put {$E_{r+1}$} at -1 .5
\put {$E_r$} at 2 4
\endpicture}
\vskip .7truecm
\centerline{\smc Figure 10}
\vskip 1truecm
\noindent
{\bf 5.5.} We introduce for $i = 1, \cdots , r$ the notation $K_i := \sum^{\kappa_i-1}_{j=0} L^{ja_i} = \frac{L^{\kappa_ia_i} - 1}{L^{a_i} - 1}$, where $-\kappa_i$ is the self--intersection number of $E_i$.  To give the general idea we first state the case $D_8$ :
\smallskip
$$\vmatrix
K_1  &  -L^{a_3}  &  L^{a_2}-1  &  0  &  0  &  0  &  0  &  0  \\
& & & & & & & \\
-L^{a_0} & K_2  &  -L^{a_1}   &  0  &  0  &  0  &  0  &  0  \\
& & & & & & & \\
0    &   -L^{a_4} &  K_3        &  -L^{a_5}  &   L^{a_4}-1  & 0 & 0 & 0 \\
& & & & & & & \\
0    &  L^{a_3}-1 &  -L^{a_2}   &  K_4  &  -L^{a_3}  &  0  &  0 &  0 \\
& & & & & & & \\
0  &  0  &  0  & -L^{a_6}  &  K_5  &   -L^{a_7}  &  L^{a_6} - 1  &  0  \\
& & & & & & & \\
0  &  0  &  0  & L^{a_5}-1  &  -L^{a_4}  &  K_6  &  -L^{a_5}  &  0 \\
& & & & & & & \\
0  &  0  &  0  &  0  &  0  &  -L^{a_8}  &  K_7  &  -L^{a_9} \\
& & & & & & & \\
0  &  0  &  0  &  0  &  0  &   L^{a_7} - 1  & -L^{a_6}  &  K_8
\endvmatrix
$$
\medskip
\noindent
In general $D_r$ is defined as the determinant of the ($r,r$)--matrix with entries $d_{i,j}$ defined as follows :
\itemitem{(1)} $d_{i,i} = K_i$ for $i = 1,\cdots,r$;
\itemitem{(2)} when $i$ is odd $d_{i,i-1} = - L^{a_{i+1}}$, $d_{i,i+1} = -L^{a_{i+2}}$, $d_{i,i+2} = L^{a_{i+1}} - 1$, and $d_{i,j} = 0$ for $j < i-1$ and $j > i+2$;
\itemitem{(3)} when $i$ is even $d_{i,i-2} = L^{a_{i-1}} - 1$, $d_{i,i-1} = -L^{a_{i-2}}$, $d_{i,i+1} = -L^{a_{i-1}}$, and $d_{i,j} = 0$ for $j < i-2$ and $j > i+1$.

\noindent
(Of course in (2) and (3) we only define $d_{i,j}$ when $1 \leq i, j \leq r$.)  When $E_0$ and
$E_{r+1}$ do not occur we replace $a_0$ and $a_{r+1}$ by 1.
\bigskip
\noindent
{\sl Remark.}  The matrix defining $D_{r-1}$ is obtained from the matrix defining $D_r$ by deleting the $r$th row and column.
\bigskip
\noindent
\demo
{5.6. About the proof of 5.4}  This can be considered as the degenerate case of the zero function $f$ in [V2, 5.6], where we studied a motivic zeta function associated to a regular function $f$ on a normal surface germ.  And the proof of [V2, 5.6] is analogous to the proof of a similar result for Igusa's local zeta functions in [V1, Sections 5 and 6].

Here a priori we just have to be careful when some discrepancy $a_i, 1 \leq i \leq r$, is zero.  (This can happen for at most one such  $i$.)  Now our Definition 3.3 is such that this does not change anything. \qed \enddemo
\bigskip
\demo{5.7. Easy example}  When $r=1$ and $a_1 \ne 0$ the contribution of $E_1$ to ${\Cal E}(S)$ is
$$\frac{(L-1)^2}{L^{a_1}-1} \big(1 + \frac{1}{L^{a_0}-1} + \frac{1}{L^{a_2}-1}\big) = \frac{(L-1)^2(L^{a_0+a_2}-1)}{(L^{a_1}-1)(L^{a_0}-1)(L^{a_2}-1)}\, ,$$
and since $a_0 + a_2 = \kappa_1a_1$ this expression is as stated in Proposition  5.4; indeed $D_1 = K_1$.  And when $a_1 = 0$ the contribution of $E_1$ is by definition $\frac{\kappa_1(L-1)^2}{(L^{a_0}-1)(L^{a_2}-1)}$, which again is as stated because now $K_1 = \kappa_1$.
\enddemo
\bigskip
\noindent
{\bf 5.8. Addendum.}  One can of course specialize Proposition 5.4 to $E(S)$ and $e(S)$; the contribution of $\cup^r_{i=1} E_i$ to $e(S)$ is $\frac{d_r}{a_0a_{r+1}}$, where $d_r$ is the absolute value of the determinant of the intersection matrix of $E_1, \cdots , E_r$.
\bigskip
\demo{5.9. Remark}  More conceptually the contributions in Proposition 5.4 and Addendum 5.8 can be considered as the essential parts of formulae for our stringy invariants of $S$ in terms of the relative log minimal model or relative log canonical model of $S$, see [V2, Sections 4 and 5]. \enddemo
\bigskip
\proclaim
{5.10. Theorem}  Let $E_i, 0 \leq i \leq r+1,$ and $D_r$ be as in 5.4 and 5.5.  We consider $D_r$ as a Laurent polynomial in $L$ (over $\Bbb Z$) with rational powers.
\itemitem{(i)} All coefficients of $D_r$ are nonnegative.
\itemitem{(ii)} The constant term of $D_r$ is (strictly) positive.  If moreover all $a_i > 0$, then this constant term is equal to 1. \endproclaim
\medskip
\demo
{Proof} (i) We first define by induction on $r$ other Laurent polynomials $A_r$ in $L$, which are useful here.  We set
$$A_1 := \sum^{\kappa_1-2}_{j=0} L^{ja_1}$$
and
$$A_r := (\sum^{\kappa_r-3}_{j=0} L^{ja_r}) D_{r-1} + L^{(\kappa_r-2)a_r} A_{r-1} \qquad \text { for } r \geq 2 ,$$
where the sum $\sum^{\kappa_r-3}_{j=0} L^{ja_r}$ is to be interpreted as zero if $\kappa_r = 2$.
\enddemo
\bigskip
\noindent
{\smc Claim.}  $D_r = (\sum^{\kappa_r-2}_{j=0} L^{ja_r}) D_{r-1} + L^{(\kappa_r-1)a_r} A_{r-1}$ for $r \geq 2$.
\bigskip
\noindent
We show the claim by induction on $r$.  For $r = 2$ we first remark that $a_0 + a_3 = (\kappa_1 - 1)a_1 + (\kappa_2 - 1)a_2$ which follows from the equalities $\kappa_1 a_1 = a_0 + a_2$ and $\kappa_2 a_2 = a_1 + a_3$.  For $r=2$ then indeed
$$\split D_2 := & K_2 K_1 - L^{a_0 + a_3} = (\sum^{\kappa_2-2}_{j=0} L^{ja_2} + L^{(\kappa_2-1)a_2})D_1 - L^{(\kappa_1-1)a_1 + (\kappa_2-1)a_2} \\ = & (\sum^{\kappa_2-2}_{j=0} L^{ja_2}) D_1 + L^{(\kappa_2-1)a_2}A_1. \endsplit$$
For $r \geq 3$ we use the following restatement of [V1, Lemma 5.6] in this context :
$$D_r = (\sum^{\kappa_r-2}_{j=0} L^{ja_r})D_{r-1} + L^{(\kappa_r-1)a_r-a_{r-1}}(D_{r-1} - D_{r-2}). \tag 6$$
The induction hypothesis yields
$$D_{r-1} = (\sum^{\kappa_{r-1}-2}_{j=0} L^{ja_{r-1}})D_{r-2} + L^{(\kappa_{r-1}-1)a_{r-1}} A_{r-2}. \tag 7$$
Combining (6) and (7) we obtain
$$\split D_r &  =  (\sum^{\kappa_r-2}_{j=0} L^{ja_r})D_{r-1} + L^{(\kappa_r-1)a_r} \!\left(\!(\frac{\sum^{\kappa_{r-1}-2}_{j=0} L^{ja_{r-1}}-1}{L^{a_{r-1}}})D_{r-2} + L^{(\kappa_{r-1}-2)a_{r-1}} A_{r-2}\! \right) \\
& = (\sum^{\kappa_r-2}_{j=0} L^{ja_r})D_{r-1} + L^{(\kappa_r-1)a_r} A_{r-1}, \endsplit$$
finishing the proof of the claim.

Finally we prove simultaneously that all coefficients of $D_r$ and $A_r$ are nonnegative, again by induction on $r$.  For $r=1$ this is obvious, and for $r \geq 2$ it is an immediate consequence of the induction hypothesis by the definition of $A_r$ and the claim.

(ii) Also this is now easily verified simultaneously for $D_r$ and $A_r$ by induction on $r$. \qed
\bigskip
\bigskip
\heading
6. Log terminal and weighted homogeneous singularities
\endheading
\bigskip
\noindent
{\bf 6.1. Hirzebruch--Jung singularities.}  Let $P \in S$ be an $A_{n,q}$--singularity with exceptional divisor $\cup^r_{i=1} E_i$ of its minimal resolution.  Then by Proposition 5.4 and Addendum 5.8 we have that
$${\Cal E}_P(S) = D_r \qquad \text { and } \qquad e_P(S) = n ,$$
where $D_r$ is as in 5.5.
\bigskip
\noindent
{\bf 6.2.} For the other log terminal surface singularities it is now also easy to give an explicit description of $E_P(S)$ as a polynomial in $uv$ with rational powers; more generally we will derive a formula for ${\Cal E}_P(S)$ for a large class of singularities, including all weighted homogeneous surface singularities.
\vskip 1truecm
\centerline{
\beginpicture
\setcoordinatesystem units <.5truecm,.5truecm>
\ellipticalarc axes ratio 3:1 360 degrees from -3 0 center at -6 0
\put {$n_1$} at  -6 0
\put {$(1)$} at  -10 0
\ellipticalarc axes ratio 3:1 360 degrees from 9 0 center at 6 0
\put {$n_k$} at  6 0
\put {$(k)$} at  10 0
\startrotation by .85 .5 
\ellipticalarc axes ratio 3:1 360 degrees from -3 0 center at -6 0
\stoprotation
\plot 0 0  -2.55 -1.5 /
\put {$n_2$} at  -5.1 -3
\put {$(2)$} at  -8.5 -5
\startrotation by -.5 .85 
\ellipticalarc axes ratio 3:1 360 degrees from -3 0 center at -6 0
\stoprotation
\plot 0 0  1.5 -2.55 /
\put {$n_i$} at  3 -5.1 
\put {$(i)$} at  5 -8.5 
\put {$E$} at 0 .6
\put {$\bullet$} at 0 0
\putrule from 0 0 to 3 0
\putrule from 0 0 to -3 0
\put {$\dots$} at -1 -4      
\put {$\dots$} at 5 -2.5     
\endpicture}
\vskip .7truecm
\centerline{\smc Figure 11}
\vskip 1truecm
\noindent
{\bf 6.3.}  Let $P \in S$ be a normal surface singularity with dual graph of its minimal log resolution $\pi : X \rightarrow S$ as in Figure 11.
Here, as in (1.6), an ellips means a chain of $\Bbb P^1$'s attached to $E$, which we denote as in Figure 12.
\vskip 1truecm
\centerline{
\beginpicture
\setcoordinatesystem units <.5truecm,.5truecm>
\put {$(i)$} at 13 0
\putrule from 0 0 to 5 0
\putrule from 7 0 to 10 0
\put {\dots} at 6 0
\multiput {$\bullet$} at  0 0  2 0  4 0  8 0  10 0 /
\put {$E$} at -.6 0
\put {$E_{r_i}^{(i)}$} at 2 -.8
\put {$E_{r_{i}-1}^{(i)}$} at 4 -.8
\put {$E_2^{(i)}$} at 8 -.8
\put {$E_1^{(i)}$} at 10 -.8
\endpicture}
\vskip .7truecm
\centerline{\smc Figure 12}
\vskip 1truecm
\noindent
We know that such a chain is determined by the two coprime numbers $n_i$ and $q_i$, where here we take $q_i$ as the absolute value of the determinant of the intersection matrix of $E^{(i)}_1, \cdots , E^{(i)}_{r_i-1}$.  Let also the central curve $E$ have genus $g$ and self--intersection number $-\kappa$.  As remarked in [D, Section 2.4], this class of singularities is very large since, by a theorem of Orlik and Wagreich [OW], it includes all weighted homogeneous isolated complete intersection singularities, for which the numbers $\{ g; \kappa; (n_1,q_1), \cdots , (n_k,q_k) \}$ are called the {\it Seifert invariants} of the singularity.  (See e.g. [D] and [W] for more information on these singularities.)

\bigskip
\proclaim
{6.4. Proposition} Let $P \in S$ be a normal surface singularity with dual graph of its minimal log resolution as in 6.3, which is not strictly log canonical.  We denote by $d$ the absolute value of the determinant of the total intersection matrix of $\pi^{-1}P$. Then

(i) the log discrepancy $a$ of $E$ is
$$a = \frac{2-2g-k + \sum^k_{i=1} \frac{1}{n_i}}{\kappa - \sum^k_{i=1} \frac{q_i}{n_i}} = \frac{\prod^k_{i=1} n_i}{d} (2 - 2g - k + \sum^k_{i=1} \frac{1}{n_i}),$$

(ii) $e_P(S) = \frac 1a (2 - 2g - k + \sum^k_{i=1} n_i) ,$

(iii) ${\Cal E}_P(S) = \frac{L-1}{L^a-1} ([E^\circ] + \sum^k_{i=1} D^{(i)}_{r_i}),$
\newline
where $D^{(i)}_{r_i}$ is the determinant of 5.5 associated to $E^{(i)}_1, \cdots , E^{(i)}_{r_i}$ and $E^{(i)}_{r_i+1} = E$ for each $1 \leq i \leq k$. \endproclaim
\medskip
\demo
{Proof}  (i) Denote for $1\leq i \leq k$ the log discrepancy of  $E^{(i)}_{r_i}$ by $a_i$. Then an easy computation, starting from $\kappa a = 2 - 2g + \sum^k_{i=1} (a_i - 1)$ and the equalities [V1, Lemma 2.4]
$$a_i = \frac{1}{n_i}(q_ia + 1) \qquad \text { for } \quad 1 \leq i \leq k,$$
yields the first expression.  The second one is essentially [A, 3.1.10].

(ii) and (iii).  Immediate from Proposition 5.4 and Addendum 5.8.  (Remark that $a \ne 0$ because $P \in S$ is required not to be strictly log canonical.) \qed \enddemo
\bigskip
\noindent
{\bf 6.5.  Log terminal singularities.}  We apply Proposition 6.4 to the singularities $P \in S$ in case (2) of (1.6(i)).  Say the $i$th ellips corresponds to an $A_{n_i,q_i}$--singularity.
\bigskip
\noindent
{\sl Case $(2,2,n_3)$} : $a = \frac{1}{\kappa n_3 - n_3 - q_3} = \frac{4}{d}$
and $e_P(S) = (n_3 + 3)(\kappa n_3 - n_3 - q_3) = (n_3 + 3) \frac d4$.
\bigskip
\noindent
{\sl Case $(2,3,3)$} : $a = \frac{1}{6 \kappa - 2q_2 - 2q_3 - 3} = \frac 3 d$
and $e_P(S) = 7(6 \kappa - 2q_2 - 2q_3 - 3) = 7 \frac d3$.
\bigskip
\noindent
{\sl Case $(2,3,4)$} : $a = \frac{1}{12 \kappa - 4q_2 - 3q_3 - 6} = \frac 2d$
and $e_P(S) = 8(12 \kappa - 4q_2 - 3q_3 - 6) = 8 \frac d2$.
\bigskip
\noindent
{\sl Case $(2,3,5)$} : $a = \frac{1}{30 \kappa - 10 q_2 - 6 q_3 - 15} = \frac 1d$
and $e_P(S) = 9(30 \kappa - 10 q_2 - 6q_3 - 15) = 9d$.
\bigskip
\noindent
Moreover ${\Cal E}_P(S) = \frac{L-1}{L^a - 1} (L - 2 + \sum^3_{i=1} D^{(i)}_{r_i})$.
\bigskip
\noindent
Since $a \in \frac {1} {\Bbb N}$ and by Theorem 5.10, we have that ${\Cal E}_P(S)$ can be considered as a polynomial in $L$ with rational powers and with nonnegative coefficients.

Combining this with (6.1) we obtain that $E_P(S)$ is such a polynomial in $uv$ for any log terminal singularity $P \in S$, as claimed in (5.2).  As a consequence, to any algebraic surface with at worst log terminal singularities, we can assign (nonnegative) stringy Hodge numbers as in (5.1).
\bigskip
\noindent
{\bf 6.6.} On the other hand let $P \in S$ be as in Proposition 6.4 but not log canonical, and suppose that $a \in \frac {1}{\Bbb Z}$.   Then necessarily $a < 0$, say $a = - \frac 1 m$ with $m \in \Bbb N$.  We have that
$${\Cal E}_P(S) = \frac{L-1}{L^{-1/m}-1} ([E^\circ] + \sum^k_{i=1} D^{(i)}_{r_i}) = -(\sum^m_{j=1} L^{j/m})( [E] - k + \sum^k_{i=1} D^{(i)}_{r_i}).$$
It is remarkable that now, again by Theorem 5.10, not $E_P(S)$ but $-E_P(S)$ is a polynomial in $u$ and $v$ (with rational powers), whose coefficients induce nonnegative \lq Hodge numbers\rq\ as in 5.1.  (In a provocative way one could say that passing from log terminal to non log canonical, through the black hole of strictly log canonical, switches positive to negative.)  This will not be true for all non log canonical $P \in S$ for which $E_P(S)$ is a polynomial in $u$ and $v$; it would be interesting to know when precisely this happens.
\bigskip
\demo{6.7. Remark}  If in (6.6) the log discrepancy of some curve $E^{(i)}_{r_i}$, intersecting the central curve $E$, is zero, then the condition $a \in \frac 1 {\Bbb Z}$ is satisfied.  Indeed, then e.g. by [V1, Lemma 2.4] we have that $a = - \frac 1{q_i}$.
\enddemo
\vskip 1truecm
\centerline{
\beginpicture
\setcoordinatesystem units <.5truecm,.5truecm>
\putrule from -2 0 to 2 0
\putrule from 0 0 to 0 -2 
\multiput {$\bullet$} at  0 0  2 0  -2 0  0 -2 /
\put {$E$} at  0 .7
\put {$E_1$} at -2.7 0
\put {$E_3$} at 2.7 0
\put {$E_2$} at -.5 -2.4
\endpicture}
\vskip .7truecm
\centerline{\smc Figure 13}
\vskip 1truecm
\demo{6.8. Example : Triangle singularities}  These are the normal surface singularities with dual graph of their minimal log resolution as in Figure 13,
where all four curves are rational and $\kappa = -E^2 = 1$.  (This is indeed a special case of (6.3).)  By Proposition 6.4(i) we have that $a = -1$, and hence that $a_1 = a_2 = a_3 = 0$.  So
$${\Cal E}_P(S) = -L(L - 2 + n_1 + n_2 + n_3) \, ,$$
where $n_i = -E^2_i$, and
$$-E_P(S) = (uv)^2  + (n_1 + n_2 + n_3 - 2)uv \, ,$$
and it is tempting to associate some kind of Hodge numbers $1$ and $n_1 + n_2 + n_3 - 2$ to $P \in S$.
\enddemo

\enddocument